\input amstex
\input amsppt.sty
\magnification=\magstep1
\hsize=30truecc
\vsize=22.2truecm
\baselineskip=16truept
\TagsOnRight
\nologo
\pageno=1
\topmatter

\def\Z{\Bbb Z}

\def\C{\Bbb C}
\def\l{\left}
\def\r{\right}
\def\b{\bigg}
\def\bg{\bigg}
\def\({\b(}
\def\[{\b[}
\def\){\b)}
\def\]{\b]}
\def\colon{{:}\;}

\def\t{\text}
\def\f{\frac}
\def\mo{\roman{mod}}

\def\em{\emptyset}
\def\se {\subseteq}

\def\eq{\equiv}

\def\ls{\leqslant}
\def\gs{\geqslant}
\def\al{\alpha}

\def\Proof{\noindent{\it Proof}}
\def\Remark{\medskip\noindent{\it Remark}}

\hbox {J. Algebra 293(2005), no.\,2, 506--512.}
\bigskip
\title A Local-Global Theorem on Periodic Maps\endtitle
\author Zhi-Wei Sun\endauthor
\affil Department of Mathematics and Institute of Mathematical Science
\\Nanjing University, Nanjing 210093, People's Republic of China
\\zwsun\@nju.edu.cn
\\ {\tt http://pweb.nju.edu.cn/zwsun}\endaffil
\abstract Let $\psi_1,\ldots,\psi_k$ be  maps from $\Z$ to an
additive abelian group with positive periods $n_1,\ldots,n_k$
respectively. We show that the function $\psi=\psi_1+\cdots+\psi_k$ is
constant if $\psi(x)$ equals a constant for
$|S|$ consecutive integers $x$ where $S=\{r/n_s:\,
r=0,\ldots,n_s-1;\, s=1,\ldots,k\}$; moreover,
there are periodic maps $f_0,\ldots,f_{|S|-1}:\Z\to\Z$ only depending on $S$ such that
$\psi(x)=\sum_{r=0}^{|S|-1}f_r(x)\psi(r)$ for all $x\in\Z$.
This local-global theorem
extends a previous result [Math. Res. Lett.
11(2004), 187--196], and has various applications.
\endabstract
\thanks 2000 {\it Mathematics Subject Classification}. Primary 20F99;
Secondary 05E99, 11A25, 11B25, 11B75, 20D60.
\newline\indent
The author is supported by the National Science Fund for
Distinguished Young Scholars (No. 10425103) and the Key Program of
NSF (No. 10331020) in China.
\endthanks
\endtopmatter
\document
\heading {1. Introduction}\endheading

In 1965 N. J. Fine and H. S. Wilf [FW] proved that for
two real sequences $\{\al_j\}_{j\gs0}$ and $\{\beta_j\}_{j\gs0}$
with respective periods $m,n\in\Z^+=\{1,2,3,\ldots\}$,
if $\al_j=\beta_j$ for all $0\ls j<m+n-\gcd(m,n)$
(where $\gcd(m,n)$ denotes the greatest common divisor of $m$ and $n$),
then $\al_j=\beta_j$ for every $j=0,1,2,\ldots$.
This result has applications in combinatorics of finite words,
see, e.g., [BB] and [R].

 The author's investigation on covers of $\Z$ by residue classes
led him to study periodic arithmetical maps in [S91], [S01], [S03]
and [S04]. Periodic maps from $\Z$ to the complex field $\C$
include Dirichlet characters and the characteristic function of a
residue class. In this paper we aim to establish
the following general local-global theorem
on periodic maps, which includes the Fine-Wilf result as a very particular case.

\proclaim{Theorem 1.1} Let $G$ be any abelian group written additively,
and let $\psi_1,\ldots,\psi_k$
be maps from $\Z$ to $G$ with periods $n_1,\ldots,n_k\in\Z^+$ respectively.
Set $\psi=\psi_1+\cdots+\psi_k$ and
$$S(n_1,\ldots,n_k)=\bigcup_{s=1}^k\bg\{\f r{n_s}:\, r=0,\ldots,n_s-1\bg\}.\tag1.1$$

{\rm (i)} There are periodic maps $f_0,\ldots,f_{|S(n_1,\ldots,n_k)|-1}:\Z\to\Z$
only depending on $S(n_1,\ldots,n_k)$ such that
$\psi(x)=\sum_{0\ls r<|S(n_1,\ldots,n_k)|}f_r(x)\psi(r)$ for all $x\in\Z$.
In particular, values of $\psi$ are completely determined by the set
$S(n_1,\ldots,n_k)$ and the initial values $\psi(0),\ldots,\psi(|S(n_1,\ldots,n_k)|-1)$.

{\rm (ii)} $\psi$ is constant if
 $\psi(x)$ equals a constant for $|S(n_1,\ldots,n_k)| \ (\ls n_1+\cdots+n_k-k+1)$
 consecutive integers $x$.
\endproclaim

\Remark\ 1.1. Theorem 1.1(i) is completely new, and it seems
that once $S=S(n_1,\ldots,n_k)$ is given those values of
$f_0,\ldots,f_{|S|-1}$ at $0,\ldots,N-1$ are uniquely determined,
where $N$ stands for the least common multiple of $n_1,\ldots,n_k$
which is also the smallest
common denominator of the rationals in $S$.
When $G$ is the additive group of a field whose
characteristic does not divide any of the periods
$n_1,\ldots,n_k$, Theorem 1.1(ii) was discovered
by the author (cf. [S04]) in May 2002 via a method rooted in [S95], though he was
unaware of the Fine-Wilf result at that time.

\medskip
 As for the cardinality of the set in (1.1), [S04, Remark 1.1] indicates that
 $$|S(n_1,\ldots,n_k)|
 =\sum\nolimits_{d\mid n_s\ \t{for some}\ s=1,\ldots,k}\varphi(d),\tag1.2$$
 where $\varphi$ is the well-known Euler function.
 By the inclusion-exclusion principle in combinatorics, we also have
$$|S(n_1,\ldots,n_k)|
=\sum_{\em\not=I\se\{1,\ldots,k\}}(-1)^{|I|-1}\gcd(n_s:\,s\in I),\tag1.3$$
because
$|\bigcap_{s\in I}\{r/n_s:\,r=0,\ldots,n_s-1\}|=\gcd(n_s:\,s\in I)$
whenever $\em\not=I\se\{1,\ldots,k\}$.

In the case $k=2$, Theorem 1.1 yields the following consequence
which is both stronger and more general than the Fine-Wilf result.

\proclaim{Corollary 1.1} Let $g$ and $h$ be maps from $\Z$
to an additive abelian group $G$ with positive periods $m$ and $n$
respectively. Then
 $\{g(x)-h(x):\,x\in\Z\}$ is contained in the subgroup of $G$ generated by
those $g(r)-h(r)$ with $0\ls r<m+n-\gcd(m,n)$; in particular,
$g$ and $h$ are identical if $g(r)=h(r)$ for all $r=0,\ldots,m+n-\gcd(m,n)-1$.
\endproclaim
\Proof. Since $|S(m,n)|=m+n-\gcd(m,n)$,
it suffices to apply Theorem 1.1 with $\psi_1=g$ and $\psi_2=-h$. \qed
\smallskip

Now we derive more consequences of Theorem 1.1.

\proclaim{Corollary 1.2}
Let $\{\psi_1(n)\}_{n\gs0},\ldots,\{\psi_k(n)\}_{n\gs0}$
be $k$ periodic complex-valued sequences
with periods $n_1,\ldots,n_k\in\Z^+$ respectively.
Let $\psi(n)=\sum_{s=1}^k\psi_s(n)$ for $n=0,1,\ldots$, and set
$$L=\sum_{\em\not=I\se\{1,\ldots,k\}}(-1)^{|I|-1}\gcd(n_s:\,s\in I).$$
Then the sequence $\{\psi(n)\}_{n\gs0}$ is a zero sequence if
its initial $L$ terms are zero. Also, the sequence $\{\psi(n)\}_{n\gs0}$
is an integer sequence if its initial $L$ terms are integers.
\endproclaim
\Proof. Observe that $L=|S(n_1,\ldots,n_k)|$ by (1.3). The first part
follows from Theorem 1.1(ii) in the case $G=\C$, and it
was realized by S. Cautis {\it et al.} [C] independent of the author's work in [S04].
The second part is a consequence of Theorem 1.1(i) in the case $G=\C$. \qed

\proclaim{Corollary 1.3} Let $G$ be an additive abelian group,
and let $c_1,\ldots,c_k\in G$ have
orders $n_1,\ldots,n_k\in\Z^+$ respectively. For any $P_1(x),\ldots,P_k(x)\in\Z[x]$,
the sum $P_1(x)c_1+\cdots+P_k(x)c_k$ vanishes for all $x\in\Z$
if it vanishes for $|S(n_1,\ldots,n_k)|$
consecutive integers $x$.
\endproclaim
\Proof. For each $s=1,\ldots,k$,
clearly $n_s$ is a period of the map $\psi_s(x)=P_s(x)c_s$.
So the desired result follows from Theorem 1.1(ii). \qed

\proclaim{Corollary 1.4} Let $\psi(x)=\sum_{s=1}^k\chi_s(P_s(x))$ for $x\in\Z$,
where each $\chi_s$ is a Dirichlet character mod $n_s$ and $P_1(x),\ldots,P_k(x)\in\Z[x]$.
Then all those $\psi(x)$ with $x\in\Z$ are linear combinations of
$\psi(0),\ldots,\psi(|S(n_1,\ldots,n_k)|-1)$ with integer coefficients.
\endproclaim
\Proof. Since $\psi_s(x)=\chi_s(P_s(x))$ is a complex-valued function with period $n_s$,
applying Theorem 1.1(i) we immediately get the required result. \qed

\proclaim{Corollary 1.5} Let $n_1,\ldots,n_k$ be positive integers,
and let $\psi:\Z\to\C$ be a function with
$\psi(0),\ldots,\psi(|S(n_1,\ldots,n_k)|)$ linearly independent over
the field of rational numbers.
Then $\psi$ cannot be written in the form $\sum_{s=1}^k\psi_s$ where
each $\psi_s$ is a function from $\Z$ to $\C$ with period $n_s$.
\endproclaim
\Proof. Let $S=S(n_1,\ldots,n_k)$.
Now that there are no $c_0,\ldots,c_{|S|-1}\in\Z$
such that $\psi(|S|)=\sum_{r=0}^{|S|-1}c_r\psi(r)$,
it suffices to apply Theorem 1.1(i) with $G=\C$. \qed

\proclaim{Corollary 1.6} Let $A=\{a_s(\mo\ n_s)\}_{s=1}^k$ be a
finite system of residue classes. If $w_A(x):=|\{1\ls s\ls k:\,
x\eq a_s\ (\mo\ n_s)\}|$ lies in a residue class $a(\mo\ m)$ for
$|S(n_1,\ldots,n_k)|$ consecutive integers $x$, then
$w_A(\Z)=\{w_A(x):\,x\in\Z\}\se a(\mo\ m)$. In particular, $A$
covers every integer an odd number of times if there are
$|S(n_1,\ldots,n_k)|$ consecutive integers each of which
is covered by $A$ an odd number of times.
\endproclaim
\Proof. Just apply Theorem 1.1(ii) with $G=\Z/m\Z$ and note that the
characteristic function of $a_s(\mo\ n_s)$ has period $n_s$.
\qed

\proclaim{Corollary 1.7} Let $A=\{a_s(\mo\ n_s)\}_{s=1}^k\ (k>1)$
be a finite system of residue classes whose maximal moduli
with respect to divisibility are distinct. Then, for any $a,b\in\Z$, we have
$$\gcd\l(w_A(a)+b,\ldots,w_A(a+|S(n_1,\ldots,n_k)|-1)+b\r)=1.\tag1.4$$
\endproclaim
\Proof. Denote the left hand side of (1.4) by $d$. Then
$w_A(a+r)\eq-b\ (\mo\ d)$ for all $0\ls r<|S(n_1,\ldots,n_k)|$,
and hence $w_A(\Z)\se -b(\mo\ d)$ by Corollary 1.6.
In view of [S05, Corollary 1.2],
$w_A(\Z)$ cannot be contained in any residue
class with modulus greater than one. Therefore $d=1$ and we are done. \qed

\heading {2. Proof of Theorem 1.1}\endheading

 The old technique used to handle the special case
 of Theorem 1.1(ii) mentioned in Remark 1.1
is invalid for the general case. Thus, we have to work along
a new line.

 Let $\Omega$ denote the ring of all algebraic integers. Clearly all roots of unity
belong to $\Omega$.

\proclaim{Lemma 2.1} Let $\psi(x)=\sum_{s=1}^kc_s\omega_s^x $ for $x\in\Z$, where
$c_1,\ldots,c_k\in\C$, and $\omega_1,\ldots,\omega_k$
are roots of unity. Suppose that
 $\prod_{\zeta\in\{\omega_1,\ldots,\omega_k\}}(x-\zeta)\in R[x]$ where
$R$ is a subring of $\Omega$ containing $\Z$. Then we have
$\psi=\psi(0)f_0+\cdots+\psi(l-1)f_{l-1}$, where
$l=|\{\omega_1,\ldots,\omega_k\}|$, and
$f_0,\ldots,f_{l-1}$ are suitable periodic maps from $\Z$ to $R$
only depending on the set $\{\omega_1,\ldots,\omega_k\}$.
\endproclaim
\Proof.  Let $\zeta_1,\ldots,\zeta_l$ be all the distinct roots of unity
 among $\omega_1,\ldots,\omega_k$, and write
$$P(z)=\prod_{t=1}^l(z-\zeta_t)=z^l-a_1z^{l-1}-\cdots-a_{l-1}z-a_l,$$
where
$$a_j=(-1)^{j-1}\sum_{1\ls i_1<\cdots<i_j\ls l}\zeta_{i_1}\cdots\zeta_{i_j}\in R
\quad \t{for}\ j=1,\ldots,l.$$

 Set $u_n=\sum_{t=1}^lc_t'\zeta_t^n$ for all $n\in\Z$,
where $c_t'=\sum_{1\ls s\ls k,\,\omega_s=\zeta_t}c_s$.
Clearly $u_n=\sum_{s=1}^kc_s\omega_s^n=\psi(n)$.
Also, $\{u_n\}_{n\in\Z}$
is a linear recurrence because
$$\align\sum_{j=1}^la_ju_{n-j}=&\sum_{j=1}^la_j\sum_{t=1}^lc_t'\zeta_t^{n-j}
=\sum_{t=1}^lc_t'\zeta_t^{n-l}\sum_{j=1}^la_j\zeta_t^{l-j}
\\=&\sum_{t=1}^lc_t'\zeta_t^{n-l}\l(\zeta_t^l-P(\zeta_t)\r)=u_{n}.
\endalign$$
If $u_{n-j}=f_0(n-j)u_0+\cdots+f_{l-1}(n-j)u_{l-1}$ for all $j=1,\ldots,l$ where
$f_0(n-j),\ldots,f_{l-1}(n-j)\in R$, then
$$\align u_n=&\sum_{j=1}^la_j\l(f_0(n-j)u_0+\cdots+f_{l-1}(n-j)u_{l-1}\r)
\\=&\(\sum_{j=1}^la_jf_0(n-j)\)u_0+\cdots+\(\sum_{j=1}^la_jf_{l-1}(n-j)\)u_{l-1}.
\endalign$$
Thus, by induction, for any $n=0,\ldots,N-1$ we can write
$u_n$ in the form $f_0(n)u_0+\cdots+f_{l-1}(n)u_{l-1}$, where
$N$ is the smallest positive integer with $\zeta_1^N=\cdots=\zeta_l^N=1$, and
$f_0,\ldots,f_{l-1}$ are suitable maps from $\{0,\ldots,N-1\}$ to $R$
only depending on the set
$\{\omega_1,\ldots,\omega_k\}=\{\zeta_1,\ldots,\zeta_l\}$.
(Actually those $f_r\ (0\ls r\ls l-1)$ can be constructed as follows:
$f_r(r)=1$, $f_r(n)=0$ for $0\ls n\ls l-1$ with $n\not=r$,
and $f_r(n)=\sum_{j=1}^la_jf_r(n-j)$ if $l\ls n<N$.)

If $x\in\Z$, and $x<0$ or $x\gs N$, then we define
$f_0(x),\ldots,f_{l-1}(x)$ to be $f_0(a),\ldots,f_{l-1}(a)$ respectively,
where $a$ is the least nonnegative residue of $x$ modulo $N$, thus
$$\psi(x)=u_a=f_0(a)u_0+\cdots+f_{l-1}(a)u_{l-1}
=f_0(x)\psi(0)+\cdots+f_{l-1}(x)\psi(l-1).$$

  In view of the above, we get the desired result. \qed
\smallskip

  The following lemma plays a central role in our proof of Theorem 1.1.

  \proclaim{Lemma 2.2} Let $\psi=\psi_1+\cdots+\psi_k$ where each $\psi_s\ (1\ls s\ls k)$
 is a complex-valued function on $\Z$ with period $n_s\in\Z^+$.
 Then $\psi$ can be written in the form
 $\sum_{0\ls r<|S(n_1,\ldots,n_k)|}\psi(r)f_r$, where
 $f_0,\ldots,f_{|S(n_1,\ldots,n_k)|-1}$ are suitable periodic maps
 from $\Z$ to $\Z$ only depending on $S(n_1,\ldots,n_k)$.
 \endproclaim
 \Proof. If $x\in\Z$ then
 $$\align \psi(x)=&\sum_{s=1}^k\sum\Sb 0\ls a<n_s\\n_s\mid x-a\endSb\psi_s(a)
 =\sum_{s=1}^k\f 1{n_s}\sum_{a=0}^{n_s-1}\psi_s(a)\sum_{r=0}^{n_s-1}e^{2\pi i\f r{n_s}(x-a)}
\\=&\sum_{s=1}^k\sum_{r=0}^{n_s-1}\(\f 1{n_s}
\sum_{a=0}^{n_s-1}\psi_s(a)e^{-2\pi i a\f r{n_s}}\)
 \l(e^{2\pi i\f r{n_s}}\r)^x.
 \endalign$$
 Observe that
 $$\prod\nolimits_{\theta\in S(n_1,\ldots,n_k)}\l(x-e^{2\pi i\theta}\r)
=\prod\nolimits_{d\mid n_s\ \t{for some}\ s=1,\ldots,k}\Phi_d(x)\in\Z[x],$$
where $\Phi_d(x)=\prod_{0\ls c<d,\,\gcd(c,d)=1}(x-e^{2\pi ic/d})$
is the $d$th cyclotomic polynomial. (That $\Phi_d(x)\in\Z[x]$
is well known, see, e.g., [IR, pp.\,194-195].)
Now it suffices to apply Lemma 2.1. \qed

\medskip
\Remark\ 2.1. Let $m,n_1,\ldots,n_k\in\Z^+$, and let $f_0,\ldots,f_{|S(n_1,\ldots,n_k)|-1}$
be as in Lemma 2.2. For each $s=1,\ldots,k$
let $\psi_s:\Z\to\Z$ be a map which has period $n_s$ modulo $m$
(i.e., $\psi_s(a)\eq\psi_s(b)\ (\mo\ m)$ whenever $a\eq b\ (\mo\ n_s)$).
Let $x$ be any integer. By Lemma 2.2 we have
$$\sum_{s=1}^k\psi'_s(x)
=\sum\nolimits_{0\ls r<|S(n_1,\ldots,n_k)|}f_r(x)\sum_{s=1}^k\psi'_s(r),$$
where $\psi_s'(x)=\sum_{0\ls a<n_s,\,n_s\mid x-a}\psi_s(a)$.
As $\psi'_s(x)\eq\psi_s(x)\ (\mo\ m)$ for each $s=1,\ldots,k$, this yields that
$$\psi(x)\eq\sum\nolimits_{0\ls r<|S(n_1,\ldots,n_k)|}f_r(x)\psi(r)\ \ (\mo\ m),$$
where $\psi=\psi_1+\cdots+\psi_k$.
For $a\in\Z$ let $\bar a$ denote the residue class $a(\mo\ m)$ in the ring $\Z_m=\Z/m\Z$.
Then we have
$$\sum_{s=1}^k\overline{\psi_s(x)}=\overline{\psi(x)}
=\sum\nolimits_{0\ls r<|S(n_1,\ldots,n_k)|}f_r(x)\overline{\psi(r)}.$$
This shows that Theorem 1.1(i) holds when $G$ is the cyclic group $\Z_m$.

 \bigskip
 \noindent{\it Proof of Theorem 1.1}. (i) Without any loss of generality,
 we simply let $G$ coincide with its subgroup generated by
 the finite set
 $$\{\psi_s(x):\, x=0,\ldots,n_s-1;\ s=1,\ldots,k\}.$$
 Since $G$ is finitely generated, there are
 $m_1,\ldots,m_l\in\Z^+$ and $n\in\{0,1,\ldots\}$ such that
 $G$ is isomorphic to the direct sum $\Z_{m_1}\oplus\cdots\oplus\Z_{m_l}\oplus\Z^n$.
 Let us identify $G$ with $G_1\oplus\cdots\oplus G_{l+n}$, where
 $G_t=\Z_{m_t}$ for $t=1,\ldots,l$, and $G_{l+1}=\cdots=G_{l+n}=\Z$.

Let $f_0,\ldots,f_{|S(n_1,\ldots,n_k)|-1}$
be as in Lemma 2.2, and let $x$ be any integer.
For $s=1,\ldots,k$
we write $\psi_s(x)$ in the vector form
$$\langle\psi_{s,\,1}(x),\ldots,\psi_{s,\,l+n}(x)\rangle,$$
where $\psi_{s,\,t}(x)\in G_t$ for $t=1,\ldots,l+n$.
Set $\psi^{(t)}=\sum_{s=1}^k\psi_{s,\,t}$ for $t=1,\ldots,l+n$.
Since $\psi_{s,\,t}:\Z\to G_t$ also has period $n_s$, we have
$$\psi^{(t)}(x)=\sum\nolimits_{0\ls r<|S(n_1,\ldots,n_k)|}f_r(x)\psi^{(t)}(r)$$
by Lemma 2.2 and Remark 2.1. Therefore,
$$\align\psi(x)=&\langle \psi^{(1)}(x),\ldots,\psi^{(l+n)}(x)\rangle
\\=&\sum\nolimits_{0\ls r<|S(n_1,\ldots,n_k)|}f_r(x)
\langle \psi^{(1)}(r),\ldots,\psi^{(l+n)}(r)\rangle
\\=&\sum\nolimits_{0\ls r<|S(n_1,\ldots,n_k)|}f_r(x)\psi(r).
\endalign$$
This proves the first part of Theorem 1.1.

(ii) Now suppose that
 $\psi(a+r)=c$ for all $r=0,\ldots,|S(n_1,\ldots,n_k)|-1$,
 where $a\in\Z$ and $c\in G$.
For $x\in\Z$ let $\psi^*(x)=\psi_s(a+x)$ for $1\ls s<k$,
$\psi^*_k(x)=\psi_k(a+x)-c$ and
$\psi^*(x)=\psi^*_1(x)+\cdots+\psi^*_k(x)=\psi(a+x)-c$.
By the first part of Theorem 1.1,
the range of $\psi^*$ is contained in the subgroup of $G$
generated by $\{\psi^*(r):\,0\ls r<|S(n_1,\ldots,n_k)|\}=\{0\}$.
Thus $\psi(a+x)-c=\psi^*(x)=0$ for all $x\in\Z$. This concludes our proof.
\qed

\enddocument